\documentclass[11pt,thmsa]{article}

\usepackage{amssymb,amsmath,amsthm,amsxtra,relsize}
\usepackage[pdftex]{graphicx}
\usepackage{subfigure}

\newenvironment{thm}{\begin{theorem}\rm}{\mbox{}\end{theorem}}

\newtheorem{theorem}{Theorem}
\newtheorem{lemma}{Lemma}
\newtheorem{proposition}{Proposition}
\newtheorem{corollary}{Corollary}
\newtheorem{definition}{Definition}
\newtheorem{example}{Example}
\newcounter{Tcap}

\def\bea{\begin{eqnarray}} \def\beaa{\begin{eqnarray*}}
\def\beq{\begin{equation}}
\def\be{\beq}

\def\aff{{\rm aff}\,}\def\Cov{{\rm Cov}\,}\def\ev{\mathbb{E}}
\def\onehalf={{\frac{1}{2}}\, }\def\Tr{{\rm Tr}\,} \def\Var{{\rm Var}\,}  \def\Proj{{\rm Proj}\,}
\def\eea{\end{eqnarray}} \def\eeaa{\end{eqnarray*}}
\def\eeq{\end{equation}}  
   
\def\ee{\eeq}
\def\Halmos{$\square$}

\begin{document}



\vspace*{.25in}
\begin{center}
{\LARGE {\bf Uses of Sub-sample Estimates to Reduce Errors in Stochastic Optimization Models} }\\
\vspace{0.5in} {\Large  John R. Birge \footnote{This work was supported by the University of Chicago Booth School of Business and by the Department of Energy under Award Number DE-SC0002587.
 Email: {\tt John.Birge@ChicagoBooth.edu}.} }\\



The University of Chicago Booth School of Business \\
5807 South Woodlawn  Avenue \\
Chicago, IL  60637, USA. \\

\today
\end{center}


%




%

%
%

\begin{abstract}
Optimization software enables the solution of problems with millions
of variables and associated parameters.  These parameters are,
however, often uncertain and represented with an analytical
description of the parameter's distribution or with some form of
sample.  With large numbers of such parameters, optimization of the
resulting model is often driven by mis-specifications or extreme
sample characteristics,  resulting in solutions that are far from a
true optimum. This paper describes how asymptotic convergence results may not be
useful in large-scale problems and how the optimization of problems based on sub-sample estimates may achieve improved results over models using full-sample solution estimates.  A  motivating example and numerical results from a portfolio optimization problem demonstrate the potential improvement.  A theoretical analysis also provides insight into the structure of problems where sub-sample optimization may be most beneficial.
\end{abstract}

\begin{flushleft}
\end{flushleft}

\section{Introduction}  \label{s:intro}

Advances in software, hardware, and algorithms for optimization have
led to orders-of-magnitude decreases in solution times and similar
increases in the sizes of problems amenable to solution (see, e.g.,
\cite{bixby}). The ability to solve problems of virtually any
scale appears to offer significant promise for operations research
methodology, but this development also comes with a price.  As
problem size increases, so does the opportunity for errors in the
parameters used to describe the model. Unfortunately, optimization
tends to focus on these errors,   leading to solutions that are
significantly inferior to a true optimal solution and, in some
cases, even inferior to na\"ive rule-of-thumb solutions that require
no optimization.

Recognizing uncertainty in optimization model parameters naturally
leads to stochastic programming formulations, but, as shown below,
these formulations also may have difficulties despite asymptotic
convergence results that suggest otherwise.  In certain examples,
such as portfolio optimization, analytical models might adequately capture the uncertainty effects in the model, but these models are
also prone to estimation errors that are inevitable in large-scale
models. One remedy for this dilemma is to use {\em robust
optimization} (see, e.g.,  
\cite{bentalnemirovski} ), which optimizes against the worst case of
a range of parameter choices.  This approach, however, still
requires some assumption on the parameter ranges and loses the form
of expected utility that is usually assumed for rational decision
making. In addition, it is difficult to characterize the effect of limited numbers of observations on reasonable uncertainty sets (for example, whether and how these sets should extend beyond the ranges of observations).

In this paper, we consider an expected utility framework directly with a focus on objectives with a reward and risk component (as in financial models) but which can include quite general utility functions. Our goals are to describe how widely solutions of sample-based problems can deviate in distribution from large-sample asymptotic results, and, under what conditions, these deviations might be most severe.  Particularly for cases in which the sample set is strictly limited (for example, by historical observations), we also wish to consider how the set of samples might best be used to estimate an optimal solution and, under what conditions, dividing the full sample into sub-sample optimization problems might provide more reliable estimates than using the full sample for a single optimization problem.

Throughout the paper, we assume that the set of parameters can be described by some
distribution (which may be used for Monte Carlo simulation) and that the modeler may not be aware of that distribution and must rely only on a set of (independent) realizations of the random variables.  The general model of the paper and assumptions appear in Section 2, where we also state
some of the known convergence properties for optimization models
in this framework. In Section 3, we explore the dependence of these properties on the dimension of the random parameters,
the effects of nonlinear dependence on the parameters, and the impact of constraints. In that section, we also present the use of sub-sample or batch estimates as a potential mechanism to reduce errors and improve convergence of solution estimates.  We then use simple examples that allow for analytical characterizations of confidence regions to demonstrate the benefit of using the sub-sample estimates.

In Section 4, we extend the analysis to a numerical study of  a stylized portfolio optimization problem that, in particular, explores the effect of the  constraints and of problem dimension on relative errors in using sub-sample estimates compared to a full-sample estimate and single optimization problem.  In Section 5, we present a general case to provide insight into the conditions that favor the use of sub-sample estimates.  Section 6 presents conclusions and directions for further investigation.

\section{Asymptotic Convergence and Universal Confidence Region Results}

The canonical problem that we consider is to find $x\in X\subset
\Re^n$ to minimize:
\begin{equation}\label{generalmodel}
\ev[f(x,\xi(\omega))],
\end{equation}
where $\ev$ represents mathematical expectation and $\omega$ is
associated with a probability space $(\Omega,\mathlarger{\mathlarger{\mathlarger{\mathlarger{\varsigma}}}},P)$.  The
function $f$ can be interpreted as the (negative of) the utility
resulting from action $x$ and outcome $\xi$.  We denote an optimal solution to \eqref{generalmodel} as $x^*$, where $x^*$ is a member of the set of optima, $X^*=\{x|\ev[f(x,\xi)]=\min_{x\in X} \ev[f(x,\xi)]\}$, with optimal value (assumed to be attainable), $z^*=\min_{x\in X} \ev[f(x,\xi)]$.  The set of parameters
are given by the random vector, $\xi:\Omega \mapsto \Xi\subset \Re^m$.  Our
interest is in cases when $m$ is large.

In some cases, (\ref{generalmodel}) can be solved directly using a
known (or supposed) distribution on $\xi$.  A typical example is the
mean-variance (
\cite{markowitz}) portfolio optimization to
minimize the portfolio variance subject to constraints on the
expected return or to minimize an objective that combines the expected return and a risk measure on the returns.  The first case can be expressed (with superscripts ${\ }^T$ denoting transposes of vectors) as $X=\{x|e^T x = 1,
 \bar{r}^T x=r_0\}$, where $e=(1,\ldots,1)^T$ is the $n$-dimensional all-ones vector, $r_0$ is the target return and $\bar{r}$
is the vector of expected returns; $f(x,\xi(\omega))=(r(\omega)^T
x-\bar{r}^T x)^2$ (i.e., such that $\xi(\omega)=r(\omega)$).  The analytical representation is then that
$\ev[f(x,\xi(\omega))]=x^T \Sigma x$, where $\Sigma$ is the variance-covariance
matrix of the returns $r(\omega)$.

The second mean-variance case of a combined objective $\ev[f(x,\xi(\omega))]=-\bar{r}^Tx+\frac{\gamma}{2}x^T \Sigma x$ has the same optimal solution form when $r(\omega)$ is normally distributed as in the case of an exponential utility such that $f(x,\xi(\omega))=\exp(-\lambda r(\omega)^T x)$ (see, for example, \cite{rjf1956}). In our simulations in Section 4, we use the mean-variance objective directly instead of the exponential form, but note the equivalence of the two objectives (up to the scaling of the risk coefficient).

While the mean-variance problem can be solved analytically if both
the mean returns $\bar{r}$ and variance-covariance matrix $\Sigma$ are
known, in practice, these distribution parameters are not known with
certainty.  In this case, if  the number of assets is $n$ , $\Sigma$ then
includes $(n(n+1))/2$ distinct elements.  Errors in any of these
estimates can lead to significant deviations in the optimum from the
(unknown) true value (see \cite{cz1993} and 
\cite{kanzhou}). In these cases, optimization may lead to worse solutions than na\"ive allocations as shown in  
\cite{demigueletal}, which found that a simplistic allocation of $x=e/n$,
equal allocation to each asset, out-performed every considered
optimization procedure for a sample set of empirical and simulated data. Even
in relatively small portfolio problems then, limited data sets can lead to
significantly sub-optimal solutions.

In the mean-variance case, the difficulty concerns  mis-estimates of
distribution parameters.  Deviations due to inexact integration or
sampling are, however, avoided.   Even when distributions are known
exactly, but finding expectations is difficult, the results can be
equally troublesome.  These difficulties occur despite encouraging
asymptotic results.

With a Monte Carlo estimated objective, we assume that the distribution of
$\xi$ is known and  consider a sample $\{ \xi^i \},i=1,\ldots,\nu$
of independent and identically distributed observations of $\xi$ that lead to the following
sample problem: \begin{equation}\label{montecarlo} \min_{x\in X}
\ \left(\frac{1}{\nu}\right)\sum_{i=1}^\nu f(x,\xi^i),
\end{equation}
which is also known as the \emph{sample-average approximation} (SAA) problem (see, e.g.,  
\cite{shapiro2003} for a general discussion).

Following the description in 
\cite{birgelouv2011},
let $x^\nu$ be the random vector of solutions to (\ref{montecarlo})
with independent  random samples.  As shown in 
\cite{kingrock} (Theorem 3.2),  the following asymptotic result suggests that sample average approximation should yield good results.

\begin{thm}\label{th1}
 Suppose that $f(\cdot,\xi)$ is convex and
twice continuously differentiable , $X=\{x|Ax\le b\}$ is a convex polyhedron,
$\nabla  f:\Re^n \times\Xi \mapsto \Re^n$:
\begin{enumerate}
\item is measurable for all $x\in X$;\\
\item satisfies the {Lipschitz} condition that there exists some
$a:\Xi\mapsto \Re$, $\int_\Xi |a(\xi)|^2 P(d\xi) < \infty$, $|\nabla
f(x_1,\xi)-\nabla f(x_2,\xi)|\le a(\xi)|x_1-x_2|$, for all
$x_1,x_2\in X$;\\
\item satisfies that there exists $x\in X$ such that $\int_\Xi |f(x,\xi)|^2
P(d\xi)< \infty$;
and, for\ $H^*=\int \nabla^2f(x^*,\xi) P(d\xi)$,\\
\item $(x_1-x_2)^T H^*(x_1-x_2)> 0, \forall x_1\not= x_2, x_1,
x_2\in X;$ \end{enumerate}
 then, for the solution $x^\nu$ to (\ref{montecarlo}), $\sqrt{\nu}(x^\nu-x^*)$ converges in distribution to $u^*$: \begin{equation}\label{clt} \sqrt{\nu} (x^\nu
-x^*) \mapsto u^*, \end{equation} where $u^*$ is the solution to:
\begin{equation}
\begin{matrix}\label{opt1}
 \mbox{min }  & \frac{1}{2} u^T H^* u + c^T u \\
\mbox{s.\ t.\ } & A_{{i\cdot}}u_{i} \le 0, i\in
I(x^*),u^{T}\nabla\bar f^{*}=0,
\end{matrix}
\end{equation}
$(x^{*},\pi^{*})$ solve $\nabla\int_{\Xi}
 f(x^{*},\xi) P(d\xi)+(\pi^{*})^{T}A=0$, $\pi^{*}\ge 0$, $Ax^{*}
 \le b$,
 $I(x^*)=\{i | A_{i\cdot}x^* = b_i\}$,
$\nabla\bar f^{*}=\int \nabla f(x^*,\xi) P(d\xi)$, and $c$ is
distributed normally $\mathcal{{N}}(0,\Sigma^*)$ with $\Sigma^*=\int (\nabla
f(x^*,\xi) -\nabla\bar f^{*})(\nabla f(x^*,\xi) - \nabla\bar f^{*})^T P(d\xi)$.
\end{thm}

This theorem implies that, asymptotically, the sample-average problem (\ref{montecarlo}) approaches a true optimal
solution to (\ref{generalmodel}) quickly with normal error distributions.  As shown, for example, in
\cite{birgelouv2011}, the convergence is actually often directly to a
point. As 
\cite{shapirohomem} discuss, such sample average approximation problems can even achieve exact convergence to an optimal solution in a finite number of samples in some cases.
These results do not, however, give an iteration number $\nu$ at
which this asymptotic regime begins to apply.  In fact, as the empirical portfolio results suggest, this regime
may only take hold for very large $\nu$.

Obtaining results that hold for any number of samples $\nu$ requires error bounds that hold in general conditions or {\it universal confidence sets} as described, for example, by 
\cite{pflug2003} and explored further by 
 \cite{vogel2008}.  A general result of this type is the following that appears in 
 \cite{daichenjrb2000}, Theorems 3.1 and 3.2 (in which $\| v\|$ refers to any consistent norm of $v\in \Re^n$).

\begin{theorem}\label{thm2} Assume that:
there exist  $a>0$, $\theta_0>0$, $\eta(\cdot): \Re^n\to \Re^+$ such that
$|f(x, \xi )|\le a\eta(\xi)$ and $\ev[e^{\theta\eta(\xi)}]<\infty$, for all $x\in X$ and $0\le \theta\le \theta_0$, then, for any $\epsilon > 0$ and for all $\nu\ge 1$, there exist $\alpha_0 > 0$, $\beta_0 > 0$, $\alpha_1>0$, and $\beta_1>0$ such that
$$P\{ |\ev_\xi[f(x^\nu,\xi)-f(x^*,\xi)]| \ge \epsilon  \}\le \alpha_1e^{-\beta_1\nu}.\label{objresult}$$
and, if $x^*$ is unique,
$$P\{\|x^\nu-x^*\|\ge \epsilon  \}\le \alpha_0e^{-\beta_0\nu}. \label{solresult}$$
\end{theorem}

This result indicates that the log of the probability of error relative to the optimal objective value or solution beyond any level is eventually linear in the number of samples. The issue in practice is how quickly this asymptotic linearity appears (or how large $\nu$ must be to overcome the effect of the constants $\alpha_0$ and $\alpha_1$ to provide meaningful bounds). 
\cite{daichenjrb2000} provides explicit results in this direction for quadratic functions.  In this paper, we will explore this convergence for specific cases and with varying tightness in the constraints.

While not discussed here, other general results to obtain confidence sets are possible with certain additional assumptions on the structure of the problem. For example, discrete decision variables may allow for certain types of bounds (see, e.g.,  
\cite{kshdm_saa_siopt2001}). Others use the relationship between the cost of an action without prior knowledge of the random realization $\xi$ compared to the {\it recourse} cost after observing $\xi$ (e.g., 
\cite{ccp2005} and 
\cite{shmoyswamysaa2004}).  In the following, we assume only general properties of $f(x,\xi)$ and the region $X$.

The focus of this paper is related to results concerning high-dimensional statistics (e.g., \cite{mjw2019}), which focus on deriving bounds that are relevant for (relatively) small numbers of samples and that often employ lower dimensional representation (such as factors for explain asset prices, e.g., \cite{sksn2023}) or effective dimension (e.g., \cite{jlpv2022}).  Another related area of interest concerns methods for reducing the loss in an objective such as in \eqref{generalmodel} when using a surrogate optimization model. A general decision analytic Bayesian view of this issue appears in \cite{jsrw2006}. Procedures for correcting the objective for both estimation and objective losses appear in \cite{aepg2022} and specifically for removing bias in the objective for a linear objective in \cite{siayrf2018}, \cite{vgpr2021}, and \cite{vgmhpr2023}, for a lasso objective in \cite{ajam2018}, and for mean-variance portfolios in \cite{kanzhou}.  Other papers that consider approaches for reducing errors include forms of cross-validation (e.g., \cite{mstone1974}) and various forms of shrinkage estimators such as the James-Stein estimator (\cite{wjcs1961}) and related forms for pooling data (e.g., \cite{vgnk2021}).

The results in this paper are also closely related to uses of bootstrapping or bagging (bootstrap aggregation) procedures in stochastic programs, such as in \cite{eichhornroemisch2007} for discrete model, and articulated for general models in \cite{lamqian2022} and explored experimentally in \cite{cw2024}. While the focus of these papers is on constructing (and evaluating) confidence intervals based on asymptotic results, the results in this paper primarily relate to performance in high-dimensional regimes where asymptotic results may have less relevance. In particular, this paper focuses in particular on the uses of independent estimates for the decision variables based on sub-sampling and their advantages for nonlinear convex objective functions $f$. The results may, however, be combined with those of the previously mentioned papers to obtain confidence intervals when relevant.

\section{Examples of Convergence Behavior and the Use of Sub-samples}

For the examples in this paper, we consider objectives that are formed from combinations of reward and risk as in the mean-variance portfolio example. As a prototypical example, we consider  a mean-risk objective of the following form:
\be\label{proto}
\min_{x\in X} \ev[-r(x,\omega)+\gamma R(x,\omega)],
\ee
where $r:X\times\Omega \to \Re$ is a \emph{reward} function (e.g.,  investment return), $R:X\times\Omega\to \Re$ is a \emph{risk} function (e.g., variance of return), and $\gamma>0$ is the \emph{risk-aversion parameter}, that represents the tradeoff between risk and return. We start with an example where the parameters determining $R$ are known and then assume progressively less information is available.

With uncertain parameters only in the reward term $r$, let $R(x,\omega)=\|x\|_1$, $X=[-1,1]^n$, and $r(x,\omega)=\xi(\omega)^T x$ to obtain the following (where we suppress the $\omega$ dependence and consider $\xi$ as a random vector):
\begin{equation}\label{examplel1}
\min_{x\in [-1,1]^n}  \ev[ -\xi^T x + \gamma \| x \|_1 ],
\end{equation}
where note that the 1-norm makes this formulation
equivalent to a linear program.

If $\ev[\xi]=0$, the optimal solution to (\ref{examplel1}) is $x^*=0$.
We would like to know the rate at which $x^\nu \to x^*$ and how that rate depends on the dimension $n$.
While the result in Theorem \ref{th1} does not apply in this case (since the objective is not twice continuously differentiable), Theorem \ref{thm2} is valid.  In fact, we can also obtain an asymptotic result directly for the sample average solution of (\ref{examplel1}), which converges  to a degenerate distribution at $x^*=0$,  a somewhat better result than the asymptotic normal distribution for the case of a quadratic risk function (e.g., $R(x,\omega)=\|x\|_2^2$) \footnote{In the case of the squared Euclidean norm risk criterion ($\|x\|_2^2$), Theorem \ref{th1} applies. The asymptotic normal distribution of $(x^\nu-x^*)$ is attained almost immediately (i.e., it only differs by the truncation of tails due to the constraints).  The probabilities of error norms greater than 1 in solutions or objective values  are the same as those given here for the $1$-norm case. Since the distribution of the solution errors can be characterized in this case, we can obtain results on the rate of convergence (or the size of the parameters in Theorem \ref{thm2}).  }.
The convergence question we wish to answer  is: when  is the sample size large enough that $\log (P\{\|x^\nu-x^*\|\ge \epsilon\})$ decreases linearly in $\nu$?

\begin{figure}[ht]
\begin{center}
\includegraphics[width=5in]{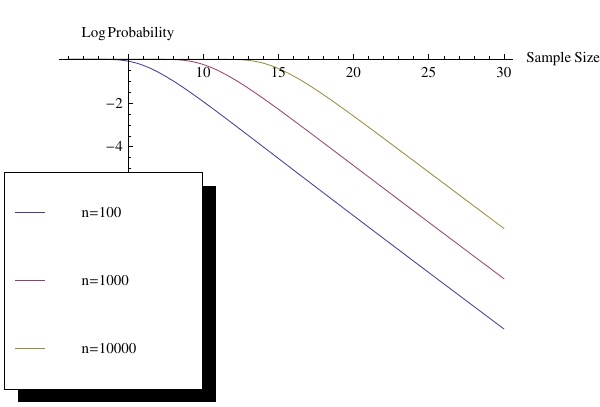}
\end{center}
\caption{Log of probability of error above one for $n=100$, 1000, and 10,000.}
\label{prob}
\end{figure}

We consider a set of samples $\{\xi^1,\ldots,\xi^\nu\}$ of the random vector $\xi$ such that each sample $\xi^i=(\xi^i_1,\ldots,\xi^i_m)$ and let $\xi^\nu=\sum_{i=1}^\nu \xi^i/\nu$.  For the sample average version of (\ref{examplel1}),   an optimal solution occurs with $\| x^\nu -
x^*\|_\infty=1$ (i.e., $x^\nu_j=\pm 1$ for some $j=1,\ldots,n$), if there exists $|\xi^\nu_j|>\gamma$ for any $j$. When
$n$ is large, the chance of $|\xi^\nu_j|\le \gamma$ for all $j$
diminishes. When each $\xi_j$ is $\mathcal{N}(0,1)$, $\xi^\nu$ is then $\mathcal{N}(0,\frac{1}{\nu})$ and $P(|\xi^\nu_j|\le \gamma)=1-2\Phi(-\gamma\sqrt{\nu})$ for $\Phi$ denoting the standard normal random variable.  The probability that this is not true simultaneously for all $n$ independent components of $\xi^\nu$ yields
\begin{equation}\label{inequality}
P\{\|x^\nu - x^*\|_\infty\ge 1\}=1-(1-2\Phi(-\gamma\nu^{0.5}))^n,
\end{equation}
where $\Phi$ is the standard normal cumulative.
To see the implication of these results, Figure \ref{prob} shows $\log (P\{\|x^\nu - x^*\|\ge 1\})$ for $\gamma=1$ for $n=100, 1000,$ and $10000$. Note also that this figure also corresponds to $P(\ev_\xi[f(x^\nu,\xi)]\ge \ev_\xi[f(x^*,\xi)]+1)$ where we note that $P$ is a probability on the random SAA solution $\xi^\nu$.

Although the asymptotic property (in terms of the probability of error above a given value) appears quickly in this case, the figure still reveals dependence on the problem dimension, which can be worse in other cases as the examples below indicate. An approach that can reduce the dimensionality effect in this case is to divide the $\nu$
samples into $K$ (in this case, non-overlapping) batches of $\nu/K$ (appropriately rounded) samples each.  As explained in Section 5, the effect of this collection of solutions from sub-samples is to reduce the overall variance of the estimated solution, which, for a convex objective, in turn yields a reduction in the out-of-sample objective in \eqref{generalmodel}.

For this approach, we let
$\xi^{\nu/K,i}$ be the mean of batch $i=1,\ldots,K$, and then solve
(\ref{montecarlo})  to obtain a solution $x^{\nu/K,i}$ for each $i$ and an overall mean solution estimate,  $\bar{x}^{\nu,K}=(1/K)\sum_{i=1}^K x^{\nu/K,i}$.

This approach, which is common in stochastic program solution estimates, e.g., see 
\cite{makmortonwood1999},  is analogous to the \emph{batch mean} method of simulation output analysis (see 
\cite{conway1963}, 
\cite{fishman1978}, and 
\cite{lawcarson1979} for origins and 
\cite{glynnwhitt1991},  
\cite{chienetal1997}, and 
\cite{steigerwilson2001} for examples of convergence results).   Alternative approaches include {\it re-sampling} procedures, such as {\it jackknife}  and {\it bootstrapping} (see, for example, 
\cite{efron1979},  
\cite{wu1990}, and 
\cite{politisromano1994} for basic results), in which, multiple overlapping sub-samples are created, and cross-validation approaches as mentioned above.   While these approaches also may lead to convergent estimates in optimization problems (see, e.g., 
\cite{eichhornroemisch2007}), bootstrapping in particular (in which, samples are chosen with replacement) may not produce consistent (and particularly non-independent) estimates  when an optimal solution $x^*$ is on the boundary of $X$ (for example, when $E(\xi)=0$ and $x\ge 0$ in \eqref{examplel1}, which follows from the analysis in the estimation bootstrapping counterexample in 
\cite{andrews2001}). Selecting sub-samples without replacement can alleviate this problem, but doing this exhaustively for multiple subsets requires many (potentially large) optimization problem solutions.  As 
\cite{goldsmanschmeiser1997} observes for simulation analysis, the total effort in such replication methods could also be spent on models with greater numbers of samples. To be consistent in terms of computation, we compare approaches using similar amounts of computational effort and assume that computational effort in calculating $x^\nu$ is roughly proportional to the number of samples $\nu$; so, that solving for $K$ sub-sample solutions, $x^{\nu/K,i},i=1,\ldots,K$, requires approximately the same computational effort as finding $x^\nu$.

In this case, the probability of a large error in the sub-sample solution estimate, $\bar{x}^{\nu,K}$, for problem \eqref{examplel1}, is:
\begin{eqnarray}\label{batch} P\{\|\bar{x}^{\nu,K}-x^*\|_\infty\ge 1\}
&\le& P\{|{x}^{\nu,i}_j|\ge 1,\forall i=1,\ldots,K;
\mbox{for\ some\ }j\in\{1,\ldots,n\},\}\\
&=& 1-(1-(2\Phi(-\gamma(\nu/K)^{0.5}))^K)^n,
\end{eqnarray}
where the probability corresponds to one minus the probability that $|\xi^{\nu,i}_j|<\gamma$ for some $i$ in $\{1,\ldots,K\}$ for each $j=1,\ldots,n$.  Note that this is also the probability of a larger than $\gamma$ error in the objective of \eqref{examplel1}.  In general, and, in particular for instances with a unique optimizer, differences in distance to optimality and the difference to optimal objective value are proportional.

The probability in (\ref{batch}) is lower than that in (\ref{inequality}) if
\be\label{normtail}
(2\Phi(-\gamma(\nu/K)^{0.5}))^K< 2\Phi(-\gamma\nu^{0.5}).
\ee
To see (\ref{normtail}), consider the asymptotic expansion of $\Phi$ as
\bea\label{asympt}
\Phi(-x)
&=& \frac{e^{-x^2/2}}{\sqrt{2\pi}}\left(\sum_{i=0}^N(-1)^ix^{-(2i+1)}\Pi_{j=0}^{\max\{i-1,0\}}(2j+1)+R^{N+1}(-x)\right),
\eea
where $R^{N+1}(-x)$ is a remainder after $N$ terms (which is positive for $N$ odd and negative for $N$ even, see, e.g., 
\cite{abramstegun1964}).
Use (\ref{asympt}) with $\nu\ge 5$, $2\le K\le \sqrt{\nu}$, and $\gamma=1$ to simplify the expressions,
\bea\label{comp}
(2\Phi(-(\nu/K)^{0.5}))^K 
&\le & \frac{\sqrt{2}}{\sqrt{\pi}}e^{-\nu/2}\left(\frac{2}{\pi}\right)^{(K-1)/2}\left(\frac{K}{\nu}\right)^{K/2}\label{in2}\\
&< & \frac{\sqrt{2}}{\sqrt{\pi}}e^{-\nu/2}\left(\frac{1}{\nu}\right)^{1/2}\left(1-(\frac{1}{\nu})\right)\label{in5}\\
&\le &  2\Phi(-\nu^{0.5})\label{in6},
\eea
where (\ref{in2}) uses (\ref{asympt}) with $N=0$ for an upper bound, (\ref{in5}) uses that $(\frac{K}{\nu})^{K/2}\le (\frac{1}{\nu})^{1/2}$ for all $ 2\le K\le \sqrt{\nu}$ and that $(\frac{2}{\pi})^{(K-1)/2}< (1-(\frac{1}{\nu}))$ for all $K\ge 2$ and $\nu\ge 5$, and (\ref{in6}) uses (\ref{asympt}) with $N=1$ for a lower bound.

The analysis in this case suggest that $K$ should be chosen sufficiently large that there is a reduction in the overall probability of a significant deviation from optimality but also that choosing $K$ too large (greater than $\sqrt{\nu}$) may allow the errors in individual sub-samples to get too great to obtain an advantage from combining their results.  The limiting value of $K=\sqrt{\nu}$ from this example provides some guidance but generally the best choice of $K$ may depend on specific problem parameters.   Some forms of cross-validation such as using the variance of $\bar{x}^{\nu,K}$ chosen over random assignments of the $K$ sub-samples for different value of $K$ may provide guidance.  In general, the choice of $K$ remains an interesting topic for further research.


This approach of sub-dividing the sample is similar to that used in 
\cite{ccp2005}, which discusses the use of multiple samples (or sub-samples) that can be ordered (due to assumed problem structure) to choose a set that avoids the worst cases and achieves a given approximation level with desired confidence.  The results in Section 2 use more general Lipschitz conditions (Theorem \ref{th1}) or  bounded exponential tail integrals (Theorem \ref{thm2}) to obtain error bounds. We now wish to explore how the rate of attaining asymptotic characteristics for both single-sample and sub-sample problem means  varies with dimension as well as greater uncertainty in the risk measure parameters and tightness in the constraints defining $X$.

\begin{figure}[ht]
\begin{center}
\includegraphics[width=5in]{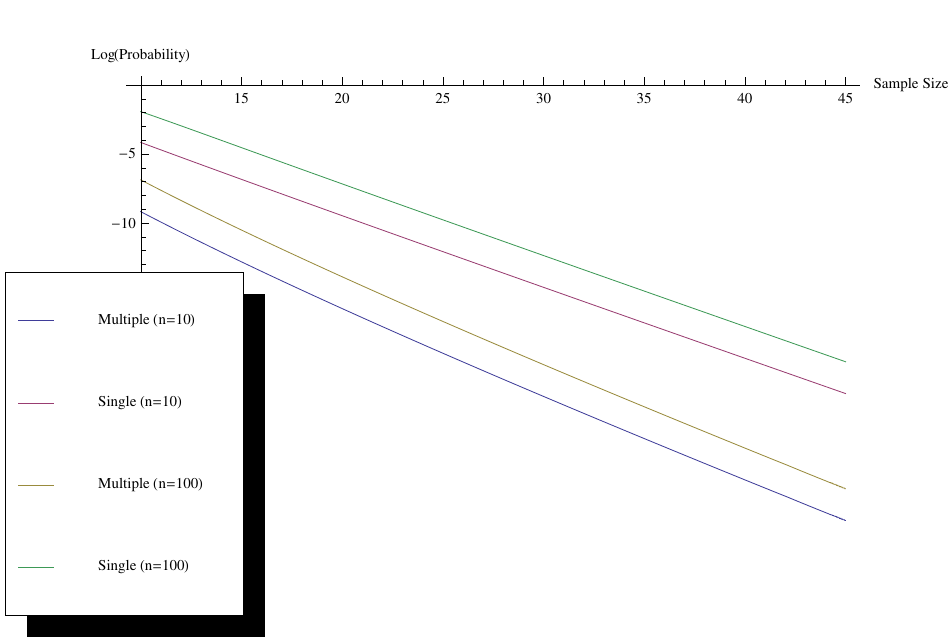}
\end{center}
\caption{Log of probability of $\infty$-norm error greater than one for $n=10$ and $n=100$ using multiple ($K=10$) and single batches.}
\label{multsingleinfnorm}
\end{figure}

To observe the effect of dimension increases, the bounds in (\ref{inequality}) and (\ref{batch}) for $K=10$ for $n=10$ and $n=100$ appear in Figure \ref{multsingleinfnorm} for the $\infty$-norm solution deviation.   Note that the error difference in dimension is approximately proportional to the increase in dimension (i.e., about two orders of magnitude in log(probability)). The multiple batch error also decreases more rapidly in sample size from an improvement relative to the single batch of approximately 5 orders of magnitude for $\nu=10$ (one sample per batch) to approximately 9 orders of magnitude for $\nu=45$.

To consider a minimal level of parameter uncertainty in the risk measure, we suppose $R(x,\omega)=\frac{x^T \Sigma x}{2}$, where $\Sigma=\sigma^2 I$, i.e., where  the risk of each asset is known to be the same and independent of all others,  but where the risk parameter, $\sigma^2$, is random\footnote{In a portfolio context, this version can also correspond to a situation where the assets are known to have identical correlation with a market portfolio but where that correlation is not known.}. Assuming an elliptical feasible region $X$ (to obtain analytical results easily), this information structure yields  the following version of the problem:
\begin{equation}\label{example2}
\min_{\|x\|_2\le 1}  \ev[ -\xi^T x + \frac{\gamma}{2} \sigma^2 \| x \|_2^2 ],
\end{equation}
where we again assume 
$\ev[\xi]=0$  and let $\ev[\sigma^2]=1$ with $\xi$ and $\sigma^2$ independent.
In the  sample average version of (\ref{example2}), the objective is:
\be\label{saa2}
-(\xi^\nu)^T x + \frac{\gamma}{2} \sigma^2_\nu \| x \|_2^2,
\ee
where
$\xi^\nu=\frac{\sum_{i=1}^\nu \xi^i}{\nu}$ and  $\sigma^2_\nu=\frac{\sum_{i=1}^\nu \sigma^2_i}{\nu}$ (where we assume that each sample observation provides an unbiased estimate of $\sigma^2$), $\| x^\nu -
x^*\|_2\ge 1$ again whenever the unconstrained  sample-average version of (\ref{example2}) has a solution, $x^{\nu,u}$,
 such that $\|x^{\nu,u}\|_2\ge 1$.    In the unconstrained case with $\gamma=1$, the asymptotic result from Theorem \ref{th1} applies, where the optimal solution $u^*=-c$, where $-c\sim \mathcal{N}(0,1)$, the standard normal distribution. This implies that
 \be\label{chi1}
 \|x^{\nu,u}-x^*\|_2^2\to  {\chi^2}(n),
 \ee
 a $\chi^2$-distributed random variable with $n$ degrees of freedom.

\begin{figure}[ht]
\begin{center}
\includegraphics[width=5in]{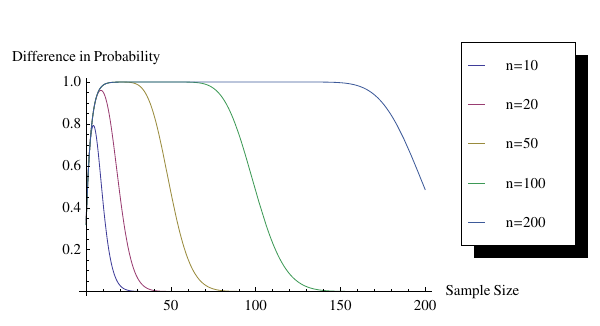}
\end{center}
\caption{Difference in error probability from the $\nu$-sample distribution and the asymptotic distribution for $n=10$, 20, 50, 100, and 200.}
\label{probdiff}
\end{figure}

 For this example, $x^{\nu,u}=\xi^\nu/\sigma^2_\nu$ so that
 \be\label{dist1}
 \|x^{\nu,u}-x^*\|_2^2=(\sum_{j=1}^n (\sum_{k=1}^\nu z_{jk}/\nu)^2)/(\sum_{i=1}^\nu y_i/\nu)^2,
 \ee
  where $z_{jk}\sim \mathcal{N}(0,1)$ and $y_i\sim \mathcal{N}(1,1)$ are independently distributed normal random variables. The result is then that
  \be\label{dist2}
  \frac{1}{\|x^{\nu,u}-x^*\|_2^2}\sim F(1,n,\nu),
  \ee
  where $F(1,n,\nu)$ is a non-central F-ratio distributed random variable with one degree of freedom and non-centrality parameter, $\lambda=E(\sum_{i=1}^\nu y_i)^2/Var(\sum_{i=1}^\nu y_i)=\nu$, in the numerator and $n$ degrees of freedom in the denominator.  From (\ref{dist2}), we have that
  \be\label{bound1}
  P(\|x^{\nu,u}-x^*\|_2\ge 1)=P( \frac{1}{\|x^{\nu,u}-x^*\|_2^2}\le 1)=P(F(1,n,\nu)\le 1).
  \ee
  Figure \ref{probdiff} shows the difference between the $\chi^2$ asymptotic probability, $P(\chi^2(n)\le 1)$, and the $\nu$-sample probability in (\ref{bound1}) as a function of $n$. The figure shows that the approach to the asymptotic result now depends more substantially on dimension than in the case with fixed risk parameter $\sigma^2\equiv 1$.  The asymptotic distribution does not yield useful confidence estimates for even moderately large sample sizes.  The number of samples necessary to obtain small confidence estimate errors from the asymptotic results increases almost linearly in dimension.

When the constraint $\|x\|_2\le 1$ is included in this example, the sub-sample approximation with $K$ batches of $\nu/K$ samples gives a solution $\bar{x}^{\nu,K}$ such that $P\{\|\bar{x}^{\nu,K}-x^*\|\ge 1\}=0$, since $\|x^{\nu/K,i}\|_2\le 1$ and $P(x^{\nu/K,i}=x^{\nu/K,k})=0$ for $i\neq k$.

\section{Portfolio Optimization Example}

The general mean-variance portfolio problem for (\ref{generalmodel}) is often written with a mean-risk objective as $f(x,\xi(\omega))=-r(\omega)^T x +\frac{\gamma}{2} (r(\omega)^Tx-\ev(r(\omega))^Tx)^2$ where $\xi(\omega)=r(\omega)$, so that, given $\ev(r(\omega))=\bar{r}$ and $\ev((r(\omega)-\bar{r})(r(\omega)-\bar{r})^T)=\Sigma$, the problem is:
 \be\label{portfolio}
 \min_{x\in X} -\bar{r}^T x + \frac{\gamma}{2} x^T \Sigma x.
 \ee

 In this example, we assume $\bar{r}$ represents expected excess returns (i.e., above a riskfree rate of return) and $x$ represents investments in $n$ different risky assets.  The entire portfolio would sum to the investor's wealth (normalized to one) (so that $1-e^T x$ would be invested in the riskfree asset).  Restricting $X=\{x|e^T x\le 1,x\ge 0\}$ would ensure a portfolio with no borrowing and no short positions in the risky assets. Allowing borrowing at the riskfree rate and no short positions would be represented by $X=\{x|x\ge 0\}$.

In general, $\bar{r}$ and $\Sigma$ are estimated from data, leading to a sample approximation version of the problem with sample estimates for both $\bar{r}$ and $\Sigma$ given as $\hat{r}=\sum_{i=1}^\nu r^i/\nu$, for $r^i$, the $i$th return observation, and $\hat{\Sigma}=\sum_{i=1}^{\nu} (r^i-\hat{r})(r^{i} -\hat{r})^T/\nu$.  Using these estimates directly to substitute for $\bar{r}$ and $\Sigma$ in (\ref{portfolio}) results in bias in the solution in the unconstrained case (see, e.g.,
\cite{kanzhou}), which can be eliminated by adding a factor $\frac{\nu}{\nu-n-2}$ in the risk term to form an adjusted problem that is unbiased in the unconstrained case as:
 \be\label{mvopt}
 \min_{x\in X} -\hat{r}^T x + \frac{\gamma \nu }{2(\nu-n-2)} x^T \hat{\Sigma} x.
 \ee

 \begin{proposition}
For $X=\Re^n$, an optimal solution, $\hat{x}$,  to (\ref{mvopt}) is un-biased, i.e., $\ev(\hat{x})=x^*$, where $x^*$ is an optimal solution of (\ref{portfolio}).
\end{proposition}

\proof
The solution to (\ref{mvopt}) is given by $\hat{x}=\frac{\nu - n - 2}{\gamma{\nu}}\hat{\Sigma}^{-1}\hat{r}$.  As noted by 
\cite{kanzhou},  $\hat{r}$  and $\hat{\Sigma}$ are independent and distributed as $\mathcal{N}(\mu, \Sigma/\nu)$ and $\mathcal{W}_n(\nu-1,\Sigma)/\nu$ respectively, where $\mathcal{N}(\mu,\Sigma/\nu)$ gives a normal distribution with mean $\mu$ and covariance matrix $\Sigma/\nu$ and $\mathcal{W}_n(\nu-1,\Sigma)/\nu$ is a Wishart distribution with $\nu-1$ degrees of freedom and covariance matrix $\Sigma$. Noting that $E((\hat{\Sigma})^{-1})=\nu \Sigma^{-1}/ (\nu-n-2)$, gives the result.\Halmos
\endproof

 When $X\neq \Re^n$, the solution of (\ref{mvopt}) may still have bias, but the objective adjustment provides more consistency in the case of interior optima.  Other adjustments of the estimated solution $\hat{x}$ can also lead to reduced errors as shown in 
 \cite{kanzhou}.  Our results are consistent with such modifications as well.

For sub-sample approximations, we use $\bar{x}^{\nu,(K)}=\sum_{i=1}^K x^{\nu/K,i}/K$ with $x^{\nu/K,i}$ that solves
\be \label{mvopt2}
\min_{x\in X} -\hat{r^i}^T x  + \gamma\frac{{ \nu/K } }{2(\nu/K- n - 2)}x^T \hat{\Sigma^i} x,
\ee
where $\hat{r^i}$ and $\hat{\Sigma^i}$ are defined analogously to $\hat{r}$ and $\hat{\Sigma}$ with samples of $\nu/K$ observations each.

While the distribution of the objective value in (\ref{mvopt}) is available (again as an $F$-distribution), analytical comparisons of the tail distributions for the objective losses and solution errors compared to optima are difficult.  Instead we illustrate the behavior with a small simulation experiment. For these results, we suppose $\nu=500$, and $K=10$ and let $\gamma=1$, $\mu=0.02e$, where $e=(1,\ldots,1)^T$, and $\Sigma=0.05*I$, where $I$ is an identity matrix.  We present the results from 1000 simulation runs for three different sets, $X$, corresponding to increasing ranges on $x$: $[0,1]^{n}$, $[-1,2]^{n}$, and $[-5,10]^{n}$. First, we consider $n=10$ and then $n=20$ to observe the effect of dimension. The results are compared relative to the optimal solution $x^*=0.4e$ in terms of relative solution distance, $\|u^\nu\|/\|x^*\|$, where $\|u^\nu\|\equiv\|{x}^\nu-x^*\|$, and relative optimal objective value for $z^*=-\bar{r}^Tx^*+\frac{1}{2}x^{*T}\Sigma x^*=-0.04$ as $(-\bar{r}^T{x}^\nu+\frac{1}{2}({x^\nu})^T\Sigma {x}^\nu-z^*)/(-z^*)= (\ev[f(x^\nu,\omega)]-\ev[f(x^*,\omega)])/(-\ev[f(x^*,\omega)])$. In the simulations, for each run of the simulation, the objective value differences represent the true expectations under the assumed distributions of the sample solutions $x^\nu$ and $\bar{x}^{\nu,K}$.

\begin{figure}[ht]
\begin{center}
\includegraphics[width=5in]{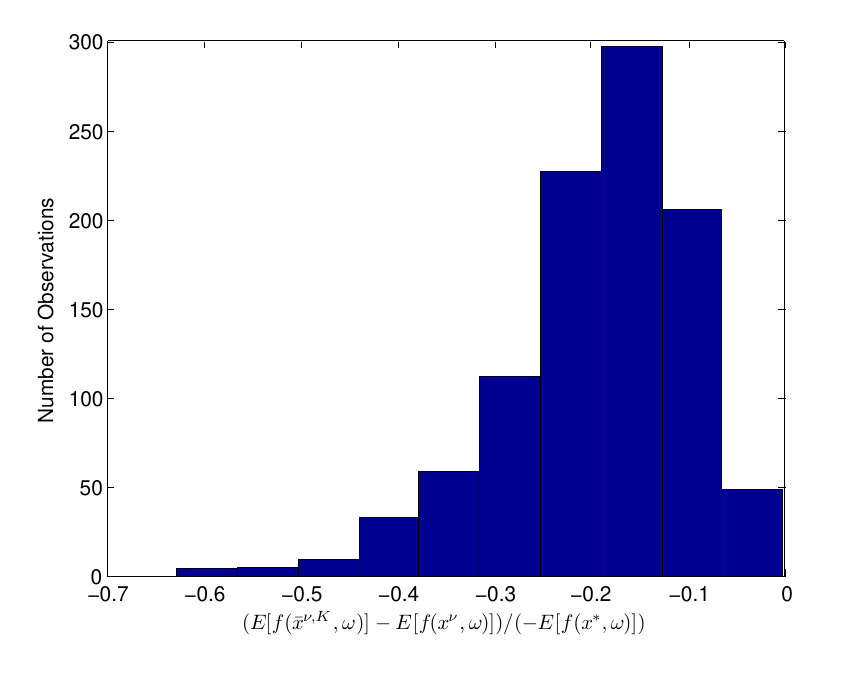}
\end{center}
\caption{Histogram of differences between relative expected objective values for sub-sample approximation minus single-sample approximation for $X=[0,1]^{n}$ and $n=10$.}
\label{hist1}
\end{figure}

\begin{figure}[ht]
\begin{center}
\includegraphics[width=5in]{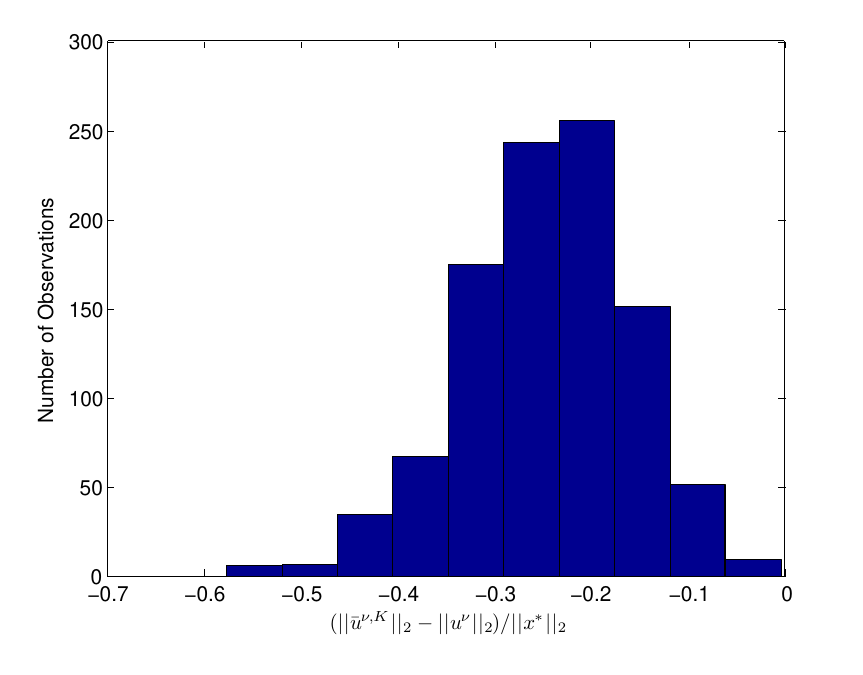}
\end{center}
\caption{Histogram of differences between relative distance from optimum for sub-sample approximation minus single-sample approximation for $X=[0,1]^{n}$ and $n=10$.}
\label{hist2}
\end{figure}

\begin{figure}[ht]
\begin{center}
\includegraphics[width=5in]{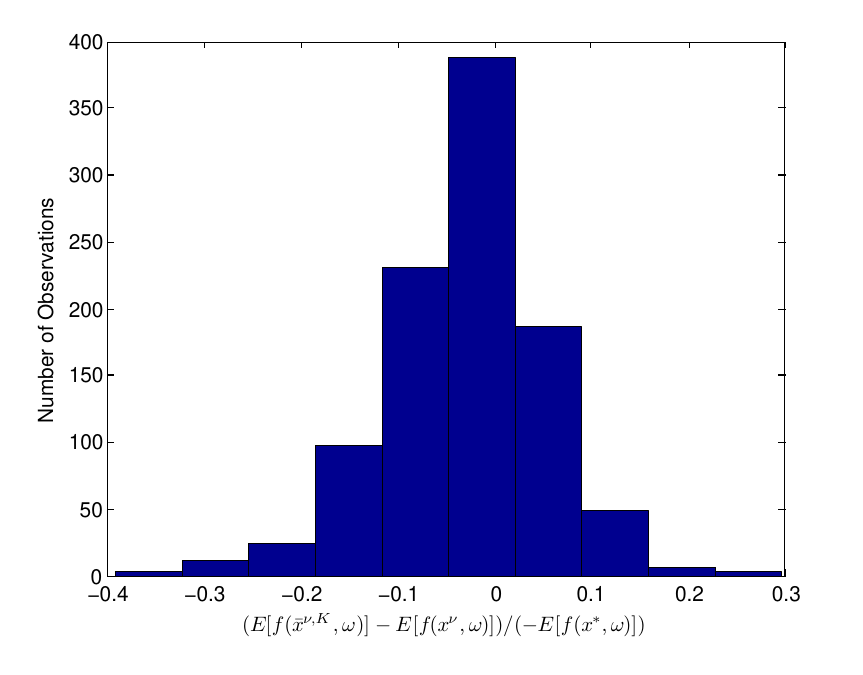}
\end{center}
\caption{Histogram of differences between relative objective values for sub-sample approximation minus single-sample approximation for $X=[-1,2]^{n}$ and $n=10$.}
\label{hist3}
\end{figure}

\begin{figure}[ht]
\begin{center}
\includegraphics[width=5in]{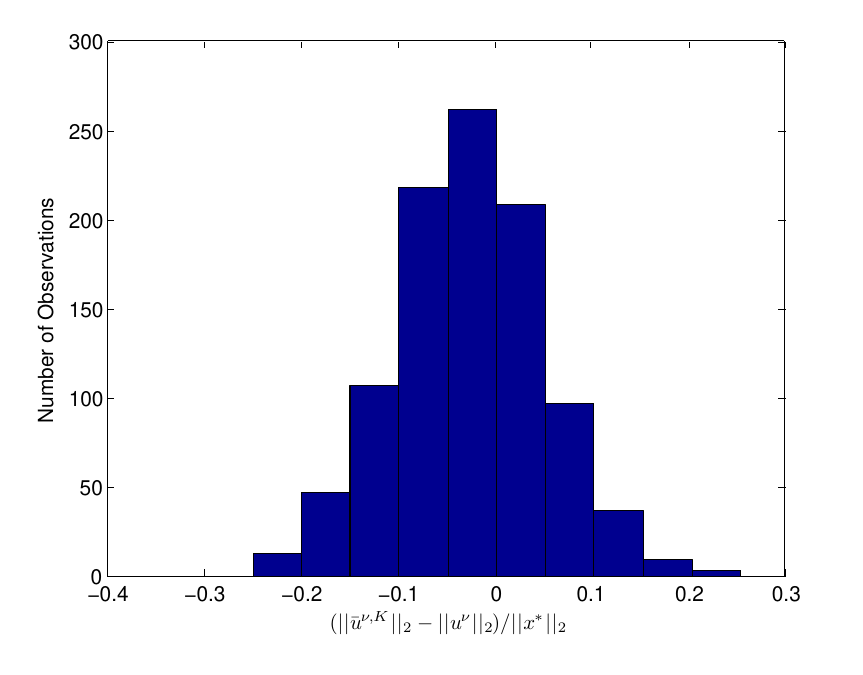}
\end{center}
\caption{Histogram of differences between relative distance from optimum for sub-sample approximation minus single-sample approximation for $X=[-1,2]^{n}$ and $n=10$.}
\label{hist4}
\end{figure}

\begin{figure}[ht]
\begin{center}
\includegraphics[width=5in]{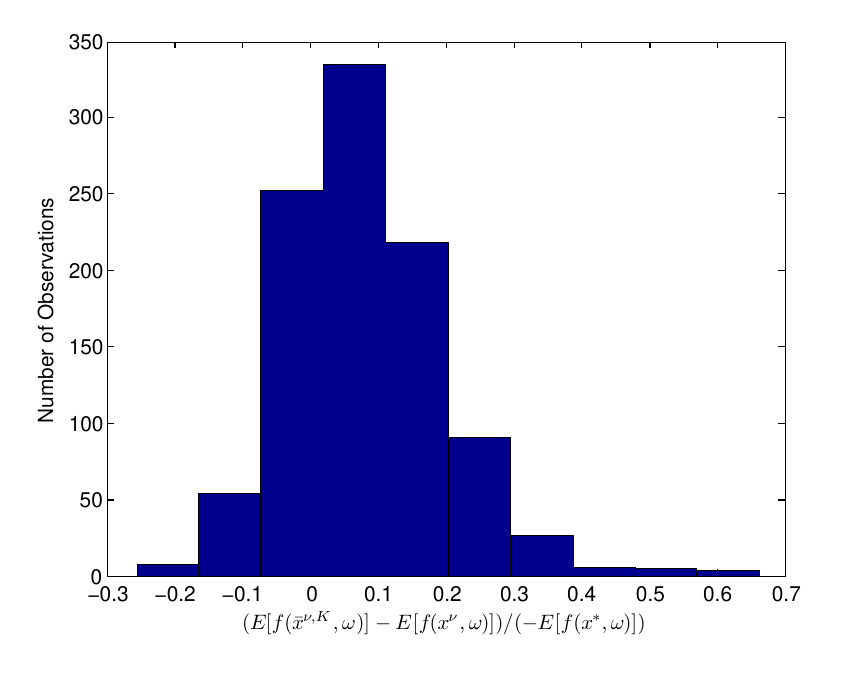}
\end{center}
\caption{Histogram of differences between relative objective values for sub-sample approximation minus single-sample approximation for $X=[-5,10]^{n}$ and $n=10$.}
\label{hist5}
\end{figure}

\begin{figure}[ht]
\begin{center}
\includegraphics[width=5in]{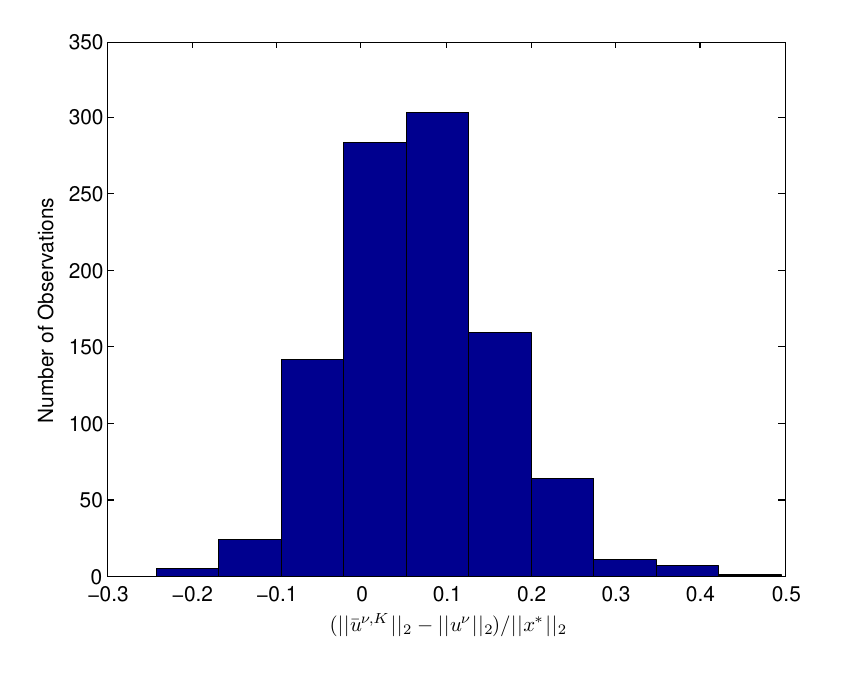}
\end{center}
\caption{Histogram of differences between relative distance from optimum for sub-sample approximation minus single-sample approximation for $X=[-5,10]^{n}$ and $n=10$.}
\label{hist6}
\end{figure}

Figures \ref{hist1}--\ref{hist6} display histograms of the results for $n=10$ for the differences in relative objective values and relative distances from the optimum between the sub-sample approximation optimal solutions and the single-sample approximation optimal solutions for the alternative values of $X$. For $X=[0,1]^{n}$ for $n=10$ in all 1000 samples, the batch approximation provided better objective values and closer approximations to the optimum.  The histogram of the differences in relative objective values appears in Figure \ref{hist1}. The histogram of the differences in relative distances from the optimum appears in Figure \ref{hist2}. While the relative objective and distance to optimality histograms are similar, the solution distance to optimality has lower mean (-25\%) compared to the difference in objective (-19\%), while the tails of the objective differences are greater than that of the distance differences to the optimal solution.  The mean of the sub-sample results in this case was also biased low with overall mean weights of $E(\hat{x}_i)=0.36$ compared to $0.40$ for the single sample overall mean weights.  This low bias caused by the constraints indicates that the adjustment to remove bias in the unconstrained case may over-compensate for bias in the constrained case (although lower weights might also be considered as compensating for uncertainty in the estimate and do have a positive overall effect as shown in the figures).

Figures \ref{hist3} and \ref{hist4} show the analogous results for the case of $X=[-1,2]^{10}$.  In this case, the greater feasible region gives a smaller advantage to the sub-sample approximation. The sub-sample approximation was better than the single sample approximation for 638 of 1000 simulation runs with average improvement relative to both the optimal value and solution of 3\%. The objective difference still has somewhat greater tails than the difference from the optimal solution. The bias in the batch approximation is reduced with average weights of 0.38.

Figures \ref{hist5} and \ref{hist6} provide the results on relative objective values and relative distance to the optimum, respectively, for $X=[-5,10]^{10}$, which is effectively an unconstrained case.  In this test, the batch approximation only improved on the single-sample approximation in 231 of 1000 runs.  Now, the single-sample approximation has an average advantage of 7\% in terms of distance to the optimum and 8\% in terms of objective value. The heavier tail observation for the objective values over the distances to optimality still appears.  In this case, in the effective absence of the constraints, the means in both the batch approximations and single-sample approximations are un-biased at $0.40$ as expected.

\begin{figure}[ht]
\begin{center} 
\mbox{\subfigure[$n=10$.]{
\includegraphics[width=3.7in]{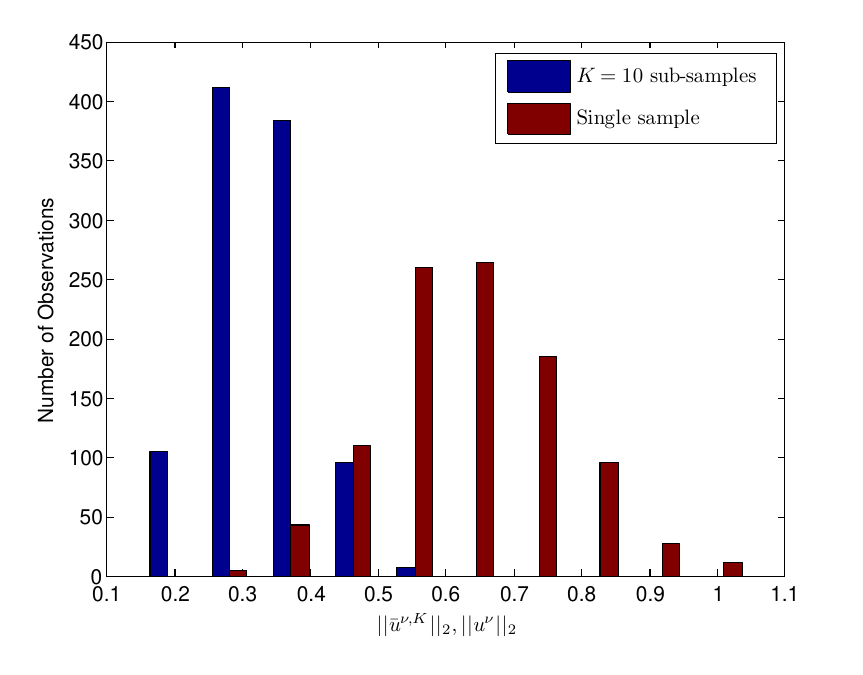}\label{hist7}}\quad
\subfigure[$n=20$.]{\includegraphics[width=3.7in]{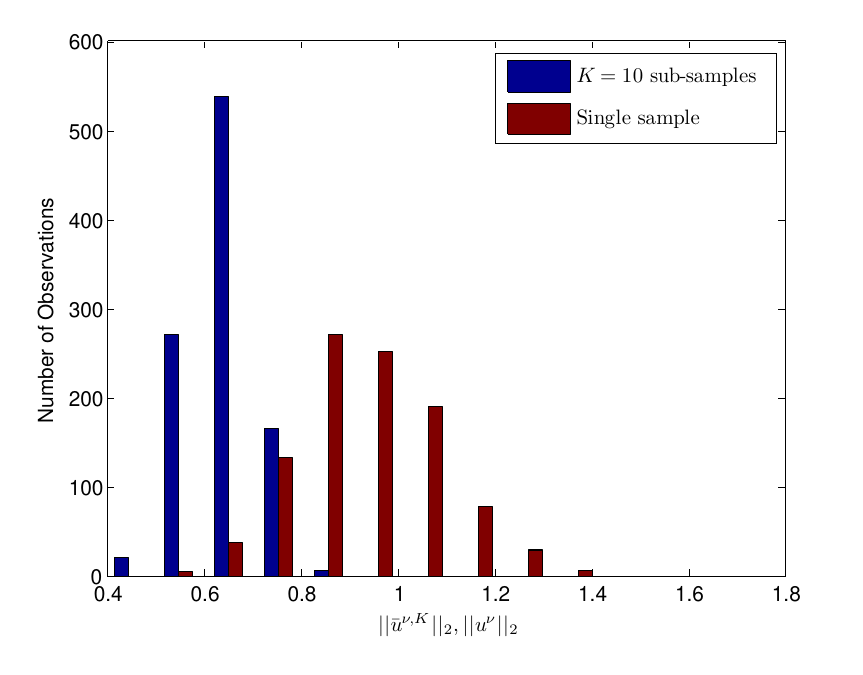}\label{hist8}}}
\end{center}
\caption{Histograms of relative distances from optimality for sub-sample versus single-sample approximations for $X=[0,1]^n$.}
\end{figure}

With $n=20$, the general results are similar to those in Figures \ref{hist1}--\ref{hist6}. For $X=[0,1]^{20}$, the sub-sample optima means were closer to optimality in 988 of 1000 runs (cf. 1000 of 1000 for $n=10$). The mean distance to the optimal solution was $17\%$ less for the sub-sample optima means than for the full-sample optima (cf. $25\%$ for $n=10$), while the difference in relative objective values was $15\%$ lower on average for the sub-sample optima mean (cf. $19\%$ for $n=10$).  To compare the sub-sample approximation results relative to those with a single-sample optimum, Figures \ref{hist7} ($n=10$) and \ref{hist8} ($n=20$) present histograms of the distances to the optimal solution for both the sub-sample optima means and the full-sample optima.  As noted, the relative differences in the distributions are similar for $n=10$ and $n=20$, but the higher dimension produces larger errors overall as expected.

\begin{figure}[ht]
\begin{center}
\mbox{\subfigure[$n=10$.]{
\includegraphics[width=3.7in]{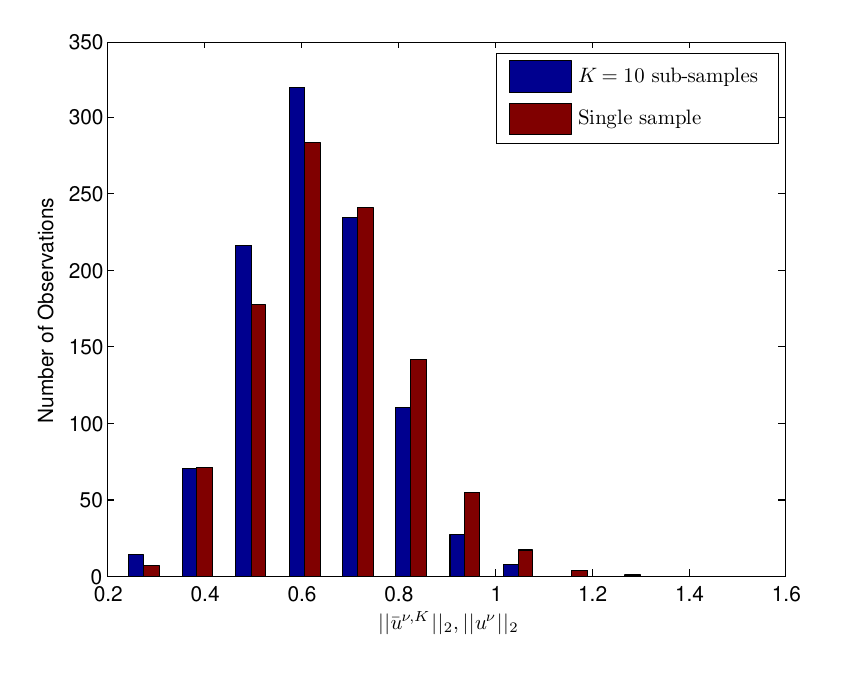}\label{hist9}}\quad
\subfigure[$n=20$.]{\includegraphics[width=3.7in]{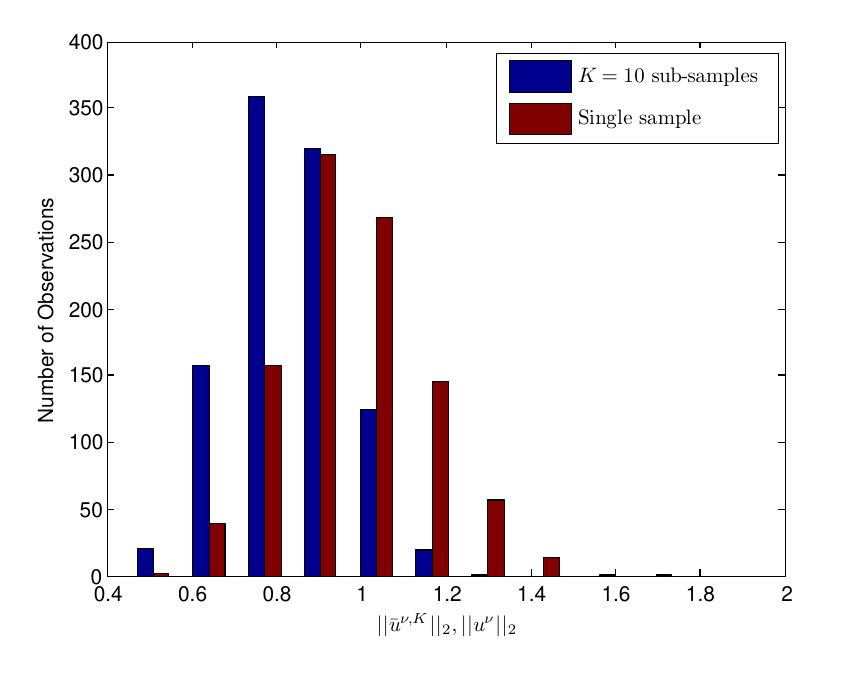}\label{hist10}}}
\end{center}
\caption{Histograms of relative distances from optimality for sub-sample versus single-sample approximations for $X=[-1,2]^n$.}
\end{figure}

For $X=[-1,2]^{n}$, a difference emerges for the $n=20$ case compared to the $n=10$ case. In this case, the sub-sample optima means are closer to optimality than the full-sample optima in 840 of 1000 runs with $n=20$ compared to 638 of 1000 runs for $n=10$. The average distance to optimality for $n=20$ is $11\%$ less for the sub-sample approximation than for the full-sample optima (cf. $3\%$ less for $n=10$) and the average relative objective loss for $n=20$ is $8.5\%$ less for the sub-sample optima (cf. $3\%$ less for $n=10$). Histograms of the distances to optimality appear in Figures \ref{hist9} ($n=10$) and \ref{hist10} ($n=20$).

\begin{figure}[ht]
\begin{center}
\mbox{\subfigure[$n=10$.]{
\includegraphics[width=3.7in]{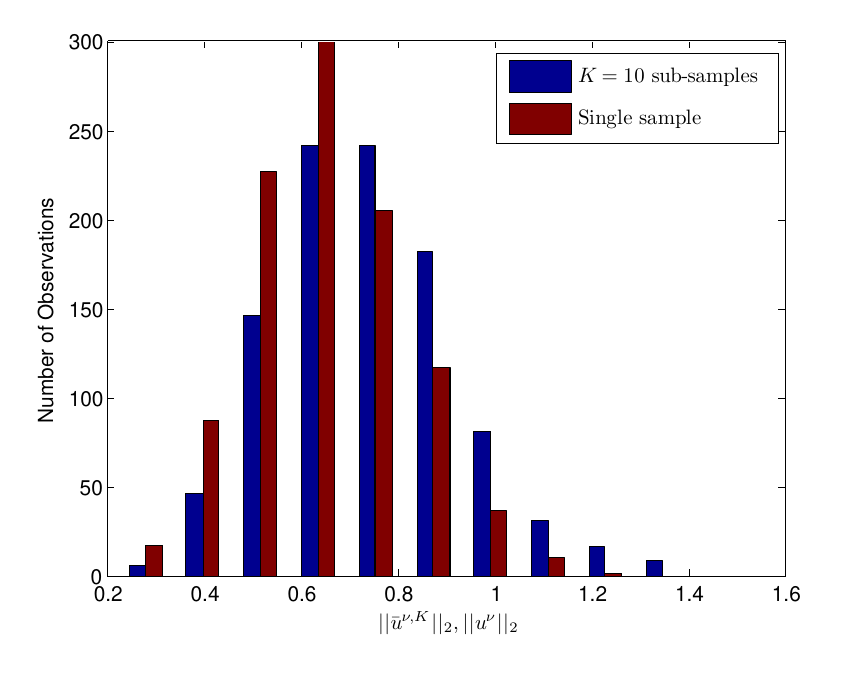}\label{hist11}}\quad
\subfigure[$n=20$.]{\includegraphics[width=3.7in]{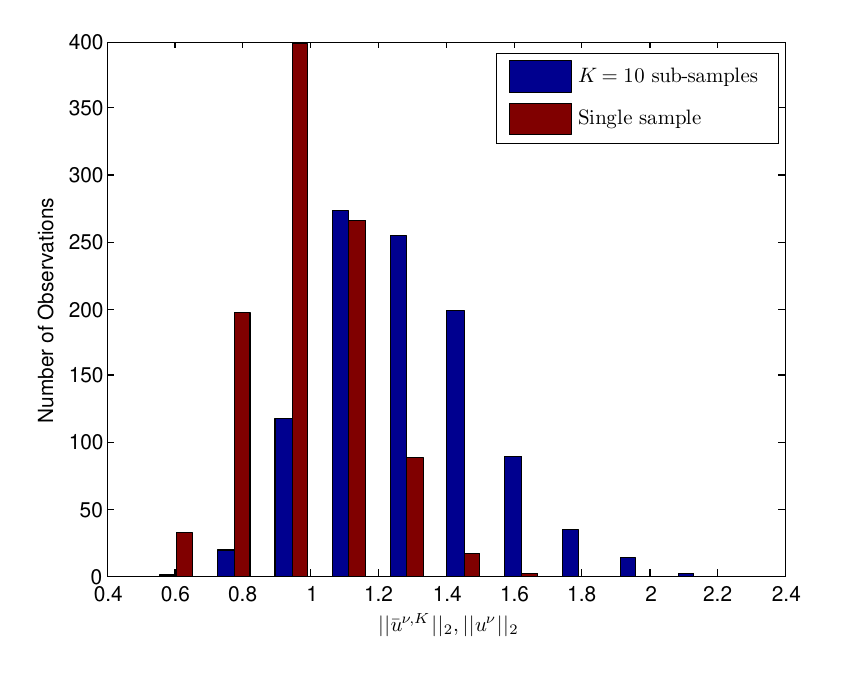}\label{hist12}}}
\end{center}
\caption{Histograms of relative distances from optimality for sub-sample versus single-sample approximations for $X=[-5,10]^n$.}
\end{figure}

For $X=[-5,10]^n$, the relative differences between the sub-sample optima means and the single sample optima becomes worse in the higher dimensional $n=20$ case than for $n=10$. For $n=20$, only 59 of 1000 runs produce lower distances to the optimum with the sub-sample approximation than with the single-sample approximation compared to 231 of 1000 runs for $n=10$. The average distance to optimality is now $17\%$ higher for the sub-sample optima means (cf. $7\%$ higher for $n=10$) and the average relative objective loss is $22\%$ higher for the sub-sample optima means (cf. $8\%$ higher for $n=10$). The histograms for these results appear in Figures \ref{hist11}--\ref{hist12}.

Overall, these results indicate  that dividing the samples into sub-sample batches and combining optimal solutions over those batches can improve results when an optimal solution should satisfy a set of constraints and may have an interior optimum (as occurred in this example). In portfolio problems, reduced errors from the addition of constraints (e.g., non-negativity for portfolio optimization) have been observed previously (e.g., 
\cite{jagma2003}). This effect is also similar to shrinking of the covariance matrix for estimation risk. While our emphasis is on the effects of constraints and sub-sample approximations for general problems, many other strategies for improving portfolio optimization estimates are possible (e.g., see the discussion in 
\cite{kanzhou}).

The results here generally confirm that errors increase in dimension for both the sub-sample and single-sample approximations and that the relative advantage of the sub-sample approximation may increase or decrease in dimension depending on the degree to which the constraints restrict the variation in sample solutions. In particular, if the constraint restrictions are very loose, then increases in dimension may not favor sub-sample approximations.  For moderate constraint restrictions, however, sub-sample approximation may improve relative to single-sample approximation as the dimension increases. Intuitively, this occurs because higher dimensional regions allow greater chances for wide variation in the single sample optimum than are possible with the averaging process of terms with limited errors in the sub-sample approximation.  The next section tries to formalize this intuition to provide further insight.

\section{General Convergence}

Error reductions using the mean of optimal solutions for sub-samples in place of a single-sample optimal solution depend on the specific characteristics of the problem instance as noted above.  We can, however, see how such improvements can become more likely as the problem dimension increases by considering the error process and the potential that errors in different sample batches do not accumulate.

We consider the canonical model in \eqref{generalmodel} and suppose a general approximation $\hat{x}^\nu(\xi^1,\ldots,\xi^\nu)$ which is formed using the observations
$\xi^{\nu\cdot}=\{ \xi^1,\ldots,\xi^\nu\}$. This approximation $\hat{x}^\nu$ (which we write without the argument  ($\xi^{\nu\cdot}$) when the context is clear) can represent the SAA problem solution or the batch-sample approximation. The expected value of the objective in \eqref{generalmodel} using $\hat{x}^\nu(\xi^{\nu\cdot})$ 
can be written (assuming the expectations below are finite) as
\begin{equation}\label{genapprox2}
\ev_{\xi^{\nu\cdot}}[\ev_\xi[f(\hat{x}^\nu(\xi^{\nu\cdot}),\xi)]]=\ev_\xi[\ev_{\xi^{\nu\cdot}}[f(\hat{x}^\nu(\xi^{\nu\cdot}),\xi)]].
\end{equation}

Assuming that $f$ is twice differentiable at a unique optimal point $x^*$ with expected gradient $\nabla\bar f^{*}=\int \nabla f(x^*,\xi) P(d\xi)$, expected Hessian $H^*=\int \nabla^2 f(x^*,\xi) P(d\xi)$, that $\ev[\hat{x}^\nu]=x^*$ (i.e., that $\hat{x}^\nu$ is an unbiased estimator of $x^*$) and $V^\nu=\Cov (x^\nu)$ is the variance-covariance matrix of the approximate solutions $x^\nu$, the expected value of $\hat{x}^\nu$ is approximately given by:

\begin{eqnarray}\label{genapprox3}
\ev_\xi[\ev_{\xi^{\nu\cdot}}[f(\hat{x}^\nu(\xi^{\nu\cdot}),\xi)]] & \approx & \ev_{\xi^{\nu\cdot}}[\ev_\xi[  f(x^*,\xi)+(\hat{x}^\nu-x^*)^T\nabla f(x^*,\xi)+\left( \frac{1}{2} \right)(\hat{x}^\nu-x^*)^T\nabla^2 f(x^*,\xi)(\hat{x}^\nu-x^*)]]\nonumber \\
&=& z^*+\ev_{\xi^{\nu\cdot}}[(\hat{x}^\nu-x^*)^T\nabla\bar f^{*}+\left( \frac{1}{2} \right)(\hat{x}^\nu-x^*)^T H^*(\hat{x}^\nu-x^*)]
 \\
&=& z^* + \left( \frac{1}{2} \right)\Tr V^\nu H^*.
\end{eqnarray}

The loss in expected objective value $L(\hat{x}^\nu)$ of using an unbiased approximation $\hat{x}^\nu$ is then approximately:
\begin{equation}\label{loss}
 L(\hat{x}^\nu)=\ev_\xi[\ev_{\xi^{\nu\cdot}}[f(\hat{x}^\nu(\xi^{\nu\cdot}),\xi)]]-z^*\approx \left( \frac{1}{2} \right)\Tr V^\nu H^*.
\end{equation}

From \eqref{loss}, we can observe that the expected objective loss from using an unbiased approximation is roughly proportional to the product of its covariance and the expected Hessian at optimality, $H^*$.  The sub-sampling procedure can be viewed as a process for reducing variance in the approximation, $x^\nu$, and hence, the expected objective loss.  As a simple example, suppose $x\in \Re^1$, $f(x,\xi)=\|x-\xi\|^2$, $X=[-1,1]$, $\xi$ is uniformly distributed on $\{-3,-1,1,3\}$, and the set of $\nu=2$ observations are two independent samples of $\xi$.  The optimal value of \eqref{generalmodel} is given by $x^*=0$ with expectation $z^*=5$.

The solutions for the standard SAA solution, $\hat{x}^\nu=x^\nu$, and the sub-sample approximation, $\hat{x}^\nu={\bar{x}}^{\nu,K}$ with $K=2$ are given in Table~\ref{ex3}.  The expected objective value $\ev_\xi[f(x,\xi)]=6$ for $x=1$ and $x=-1$; thus, the expected loss from the SAA solution is $L(x^\nu)=0.75$ while the expected loss from the sub-sample approximation is $L({\bar{x}}^{\nu,K})=0.5$. These losses coincide with the variances of these solutions: $\Var(x^\nu)=0.75$ and $\Var({\bar{x}}^{\nu,K})=0.5$.

\begin{table}
\begin{center}
\begin{tabular}{|c||c|c|c|c|}
  \hline
 & \multicolumn{4}{|c|}{$\leftarrow$ $\xi^1$ $\rightarrow$}\\
    \hline
 $ \xi^2 $  ($\big\downarrow $)     & -3 & -1 & 1 & 3 \\
    \hline
 -3 & (-1,-1) & (-1,-1) & (-1,0) & (0,0) \\
 -1   & (-1,-1) & (-1,-1) & (0,0) & (1,0) \\
 1   & (-1,0) & (0,0) & (1,1) & (1,1) \\
 3  & (0,0) & (1,0) & (1,1) & (1,1) \\
  \hline
\end{tabular}
\end{center}

\caption{Solutions $(x^\nu,{\bar{x}}^{\nu,K})$ for the SAA approximations  for $\nu=2$, $K=2$  for realizations $\xi^1$ and $\xi^2$ of \eqref{generalmodel} with $f(x,\xi)=\|x-\xi\|^2$, $X=[-1,1]$, and $\xi$  uniformly distributed on $\{-3,-1,1,3\}$.}
\label{ex3}
\end{table}

This example illustrates how the variance of the averages of sub-sample solutions to stochastic optimization problems often have lower variance than the solutions of full-sample SAA problems.  In many cases, if the SAA objective has an unconstrained optimum $x^{\nu,unc}$, then the constrained solution, $x^{\nu}=\Proj_{X}( x^{\nu,u})$ for some form of projection of the unconstrained solution onto the constraints.  The sub-sample average solutions can be viewed as averaging solutions
$x^{\nu/K}=\Proj_{X}(x^{\nu/K,u})$.  While each of the unconstrained solutions $x^{\nu/K,u}$ may have higher variance than $x^{\nu}$, the means of the sub-sample optima have lower variance.  Using overlapping samples (as in cross-validation) induces correlation among the $x^{\nu/K,u}$ solutions, potentially increasing the variance depending on the relative advantage  of increasing samples to the disadvantage of  losing  independence among the sub-samples.

For this explanation, we assumed that unbiased solutions were available as potentially in the case of mean-variance portfolio optimization as discussed above.  Debiasing procedures in other cases that fit  the case of a convex objective with constraints may also be available.  For example,  lasso regression for fitting observations $(\xi^i_0,\xi^i_1,\ldots,\xi^i_n)$ in which $\xi_0$ is assumed to be a linear function of features $(\xi_1=1,\xi_2,\ldots,\xi_n)$ plus noise can be written as \eqref{generalmodel} with
$f(x,\xi)=\|\xi_0- (\xi^i_1,\ldots,\xi^i_n)x\|^2$ and $X=\{x|\|x\|_1\le \alpha\}$ for some $\alpha$.  The example of lasso regression and other forms of regularization can be combined with the sub-sample estimation as an additional form of regularization to achieve further gains over the benefit from one type of regularization alone.

Procedures for debiasing the solutions such as those described in \cite{ajam2018} can then be used to ensure low bias in the sub-sample solutions $x^{\nu/K}$.  In general,  $x^{\nu/K}$ may have bias in which cases, the sub-sample approach must also consider the relative impact of the bias of using smaller samples versus the potential from reduced variance in the approximation.  The following discussion considers this more general situation and gives some additional perspective on the advantages of using sub-samples.

Suppose that the error for the sample average problem with $\nu$ samples is $u^\nu$, where
\be\label{errordef}
u^\nu= x^\nu-x^*.
\ee
For batching $\nu$ samples into $K$ groups to be effective, we need
\be\label{improve}
\ev_\xi[f(\sum_{i=1}^K x^{\nu/K,i}/K,\xi)]=\ev_\xi[f(x^*+\sum_{i=1}^K u^{\nu/K,i}/K,\xi)]\prec \ev_\xi[f(x^*+u^\nu,\xi)]=\ev_\xi[f(x^\nu,\xi)],
\ee
where $\prec$ is used to indicate ordering in terms of some general loss function with respect to the sampling distribution.   This ordering is generally consistent with an ordering based on the magnitude of the error; so, an alternative goal is then
\be\label{improve2}
\|\sum_{i=1}^K u^{\nu/K,i}/K\|\prec \|u^\nu\|.
\ee

To see how (\ref{improve2}) can arise, consider the form of the asymptotic distribution in Theorem \ref{th1}.  For $x^*$ in the interior of a face $F^*$ of $X$ of dimension $N$, the asymptotic distribution  of the error includes optimization of a quadratic function plus the random linear function $c^T u$ where $c$ is normally distributed as in the theorem.  Under non--degeneracy, $\pi^*_{I(x^*)}>0$, so that $A_{I(x^)}u_{I(x^*)}=0$ (from substituting for $\pi^*_{I(x^*)}A_{I(x^*)}=\nabla \bar{f}^*$ in $(\nabla \bar{f}^*)^T u=0$).  The asymptotic distribution ${u}_a^\nu$ and an associated multiplier ${\pi}^\nu$  then solve
\be
\begin{pmatrix}
 H^* & A_{I(x^*)}^T\\
A_{I(x^*)} & 0
\end{pmatrix}({u}_a^\nu; {\pi}^\nu)=(-c;0).
\ee

 The distribution of ${u}_a^\nu$ is consequently distributed normally in the  subspace of dimension $N$, where $A_{I(x^)}u_{I(x^*)}=0$. For the actual error $u^\nu$, $x^\nu$ must also be feasible, which leads $u^\nu$ to be (asymptotically) a truncation of $u_a^\nu/\sqrt{\nu}$ such that $x^\nu$ is feasible. For wide variations in $c$, this leads to concentrations on the extreme points of $X$ in $F^*$.  For bounded $X$, we have at least $N+1$ such extreme points. The worst case concentration of errors then occurs when errors are concentrated in the $N+1$ directions of these points.  Each error $u^i$ from the sub-problem for  batch $i$ is an independent sample  with the truncated distribution of $u_a^{\nu/K}/\sqrt{\nu/K}$, which then corresponds to a random selection among at least $N+1$ directions. If the expected weights on these directions are roughly proportional, then averaging the batch-optimal solutions can lead to lower errors than using a single sample.

 To make this concept more formal, suppose that
 \be\label{simplex}
 F^*\subset \hat{F}=\{x|x=x^*+\sum_{i=1}^{N+1}\lambda_i d_i,\sum_{i=1}^{N+1} \lambda_i\le 1; \lambda_i\ge 0,i=1,\ldots,N+1\},
  \ee
  where $\hat{F}$, a simplex enclosing $F^*$, is determined by
  $D=\{d_i,i=1,\ldots,N+1\}$, a set of directions such that $d_i^Td_j\le 0, \forall i\neq j$ and $\|d_i\|\le M, \forall i$.

For any error realization, $ {u}^\nu={x}^\nu-x^*$, such that ${x}^\nu\in F^*$,
 \be
 \|{u}^\nu\|^2=(\sum_{i=1}^{ N+1} {\lambda_i^\nu} d_i)^T(\sum_{i=1}^{N+1}  {\lambda_i^\nu} d_i)\le \max_{1\le i\le {N+1}} \|d_i\|_2^2\le M^2.
 \ee

Now, suppose a set of $K$ samples with errors, $\{u^{\nu/K,1},\ldots,u^{\nu/K,K}\}$, and let $\bar{u}^{\nu,K}=\sum_{i=1}^K u^{\nu/K,i}/K$, 
such that $\|\bar{u}^{\nu,K}\|_2^2= $
\bea
(\bar{u}^{\nu,K})^T \bar{u}^{\nu,K} & =& (\sum_{j=1}^{N+1}(\sum_{i=1}^{K} \lambda^{\nu,K}_{ij})d_j)^T(\sum_{j=1}^{N+1}(\sum_{i=1}^K \lambda^{\nu,K}_{ij})d_j)/K^2\\
&\le & (\sum_{j=1}^{N+1}(\sum_{i=1}^K \lambda^{\nu,K}_{ij})^2 (d_j)^Td_j)/K^2.\label{normub}
\eea

Now, suppose that  $|E(\lambda^{\nu,K}_{ij})|\le  g/N$ for some constant $g>0$. We then have, for each $j=1,\ldots,N+1$,  the multiplier on the error norms in (\ref{normub}), $\sum_{i=1}^K \sum_{l=1}^K\lambda^{\nu,K}_{ij}\lambda^{\nu,K}_{lj}$, 
can be bounded to obtain an overall bound on the error tails as in the following proposition.

\begin{proposition}\label{universalbound}
Suppose $\{x^{\nu/K,i},i=1,\ldots,K\}$ are solutions of sample problem (\ref{montecarlo}) with $K$ independent sets of  $\nu/K$ independent samples each with solution bias, $b_{\nu/K}=\|\ev[x^{\nu/K,i}-x^*]\|$, where the feasible region, $X\subset\Re^n$, is compact, and $x^*$ is an optimal solution of (\ref{generalmodel}), where $x^*\in F^*$, such that $\dim (\aff F^*)=N+1$ and $F^*\subset \hat{F}$, as in (\ref{simplex}), with $d_i^T d_j\le 0, i\neq j$, and  $\|d_i\| \le M, i=1,\ldots,N+1$,  and  
${x}^{\nu/K,i}=x^*+\sum_{i=1}^{N+1} \lambda^{\nu,K}_i d_i$, $\bar{x}^{\nu,K}=\sum_{i=1}^{K}{x}^{\nu/K,i}$, and that $\ev[\lambda^{\nu,K}_i]\le \frac{g}{N}$ 
and $g>0$; then, for any $a>0$, the error, $\bar{u}^{\nu,K}=\bar{x}^{\nu,K}-x^*$, in the average of the $K$ batch sample average problems satisfies
\be\label{univ}
P\left(\|\bar{u}^{\nu,K}\|\ge b_{\nu/K}+\frac{aM((N+1)g(N-g))^{1/2}}{K^{1/2}N}\right)\le \frac{1}{a^2+1}.
\ee

\end{proposition}

\proof
From (\ref{normub}), we assume that
\be\label{order}
\|\bar{u}^{\nu,K}\|\prec \frac{(N+1)^{1/2}M(\sum_{i=1}^K \hat{\lambda}_{i})}{K},
\ee
where $\prec$ indicates stochastic ordering and $\hat{\lambda_{i}}\succeq \lambda^{\nu,K}_{ij}$ for all $j=1,\ldots,N+1$. In particular, let $\hat{\lambda_{i}}$ be a Bernoulli random variable with expectation $\frac{g}{N}$ to obtain this ordering.  Now, $\sum_{i=1}^K \hat{\lambda}_{i}$ is Binomial($\frac{g}{N},K)$. The standard deviation of the random variable on the right-hand side in (\ref{order}) is then $\frac{((N+1)g(N-g))^{1/2}M}{K^{1/2}N}$. Applying the one-sided Chebyshev inequality then yields the result in (\ref{univ}).\Halmos
\endproof

The result in Proposition \ref{universalbound} remains valid for $F^*=X$ and $n=N$, although useful bounds on $\ev[\hat{\lambda}_i]$ are more likely to depend on the lowest dimensional face containing $x^*$ as given here. As an example of specific results, for unbiased formulations ($b_{\nu/K}=0$), $g<<N$, and $a=K^{1/4}$, the error is greater than $\frac{g^{1/2}M}{K^{1/4}}$ with probability at most $\frac{1}{\sqrt{K}+1}$.     The bias, $b_{\nu/K}$, can be bounded under certain regularity conditions (see, e.g.,  
\cite{roemischschultz1991}) as $O((\nu/K)^{-\frac{1}{2n}})$ in general, and, as in the examples above, the bias can be reduced or eliminated.

Note that the bound in Proposition \ref{universalbound} can apply when the single-sample errors may not be useful.  For example, in the case of extremely high variance in $\Sigma^*$ in the asymptotic problem (\ref{opt1}), $\lambda^\nu_j$ may be approximately Bernoulli distributed when $F^*=X$ and $\|d_i\|=M, \forall i$, so that  it is possible that $P\{\|x^\nu-x^*\|\ge \delta \}\approx 1$ for any $\delta< M$, while the result in (\ref{order}) can obtain bounds that increase in precision and confidence in $K$.

\section{Conclusions}

Large-scale stochastic optimization problems present significant
problems for effective solution. The goals in this paper have been to show that problem dimension has different effects on asymptotic convergence properties for these problems depending on the structure of the objective and constraints and that combining optimal solutions from sub-samples can improve Monte Carlo sample estimates in cases where the asymptotic regime may not be relevant.

 The results here indicate that, in cases with linear effects of the random parameters on the sample problem solution,  exponential decreases in error probabilities start to appear for sample sizes that are less than linear in the dimension, but that, when random parameters have interaction effects, the appearance of such error reductions  can deteriorate rapidly in dimension.  This first set of observations suggests that modelers should exercise caution in using asymptotic results to construct confidence regions for sample average approximation problems.  For example, trying to estimate the parameters in (\ref{opt1}) to construct a confidence region from a single sample may not be warranted.

Using sub-samples, however, can provide a more definitive test of whether a consistent set of parameters applies to (\ref{opt1}) and (\ref{clt}).  Increasing sizes of these samples can also be used to estimate the parameters for the bounds in Theorem \ref{thm2}. Even if these stronger convergence properties cannot be identified, increasing sub-sample sizes may be used to test for bias and possibly establish weaker bounds such as those in (\ref{univ}).

The advantages of sub-sample estimates appears particularly in convex optimization models when low-bias approximate optima can be obtained from samples and sub-sample averages have lower variance than those from using the full sample average approximation solutions.  This observation indicates the possibility of using bootstrap estimates of these quantities to achieve overall confidence intervals on out-of-sample optimal values, as mentioned in the introduction and the discussion of results related to bootstrapping and bagging methods.   Investigating the construction of estimates based on these ideas more thoroughly as well as designing schemes for optimal choices of the batch number $K$ are promising areas for future research.

\section{Compliance with Ethical Standards}

The author attests to having no potential conflicts of interest in this manuscript.  Part of this work was developed under a research grant from the US Department of Energy and supported by the University of Chicago Booth School of Business as referenced in the acknowledgements.  The opinions expressed in this paper are those of the author and do not represent those of the US Department of Energy, the University of Chicago, the University of Chicago Booth School of Business, or any other organization.  This work has not involved human subjects or the use of animals.



\bibliographystyle{plain}

\bibliography{stoprorefs}

\end{document}